\documentclass[11pt]{amsart}
\usepackage[a4paper,margin=1.2in]{geometry}
\usepackage{amsmath,amssymb,amsthm,mathtools}
\usepackage{tikz}
\usepackage[numbers]{natbib}
\usepackage[colorlinks=true,linkcolor=blue,citecolor=blue,urlcolor=blue]{hyperref}
\usepackage[nameinlink,capitalize]{cleveref}

\numberwithin{equation}{section}

\theoremstyle{plain}
\newtheorem{theorem}{Theorem}[section]
\newtheorem{lemma}[theorem]{Lemma}
\newtheorem{proposition}[theorem]{Proposition}

\theoremstyle{definition}
\newtheorem{remark}[theorem]{Remark}

\crefname{theorem}{Theorem}{Theorems}
\crefname{lemma}{Lemma}{Lemmas}
\crefname{proposition}{Proposition}{Propositions}
\crefname{remark}{Remark}{Remarks}
\crefname{section}{Section}{Sections}
\Crefname{section}{Section}{Sections}

\newcommand{\R}{\mathbb R}
\newcommand{\E}{\mathbb E}
\newcommand{\PP}{\mathbb P}
\newcommand{\one}{\mathbf 1}

\newcommand{\bPhi}{\overline\Phi}

\newcommand{\doilink}[1]{\href{https://doi.org/#1}{doi:#1}}

\begin{document}

\title[Multivariate self-normalized moderate deviations]
{A geometric proof of multivariate self-normalized moderate deviations via ridge regularization}

\author{Jing Bing-Yi}
\address{School of Artificial Intelligence,
The Chinese University of Hong Kong, Shenzhen,
Shenzhen 518272, China}
\email{bingyijing@cuhk.edu.cn}

\author{Liu Yiming}
\address{School of Economics, Jinan University, Guangzhou 510632, China}
\email{liuyiming@jnu.edu.cn}

\author{Wang Shaochen}
\address{School of Mathematics, South China University of Technology, Guangzhou 510640, China}
\email{mascwang@scut.edu.cn}

\author{Zhou Wang}
\address{Department of Statistics and Data Science, National University of Singapore, Singapore 119077, Singapore}
\email{stazw@nus.edu.sg}



\begin{abstract}
Let $X,X_1,X_2,\ldots$ be i.i.d.\ random vectors in $\mathbb R^d$, where $d\ge1$ is fixed.  Assume that $\E X=0$, $\E XX^\top=I_d$ and $\E\|X\|^3<\infty$.  We prove  a Cram\'er type moderate deviation theorem for a multivariate self-normalized radius statistic, that is, for every $a_n=o(n^{1/6})$,
\[
   \sup_{0\le x\le a_n}
   \left|\frac{\PP(U_{n,+}\ge x)}{\PP(\chi_d\ge x)}-1\right|\to0,
\]
where $\chi_d$ is the Euclidean norm of a standard normal vector in $\R^d$ and
\[
   U_{n,+}^2=S_n^\top A_n^+S_n,
   \qquad S_n=\sum_{j=1}^nX_j,
   \qquad A_n=\sum_{j=1}^nX_jX_j^\top.
\]
Here $A_n^+$ denotes the Moore--Penrose inverse.  We further prove that the same result holds for the corresponding Hotelling $T^2$ statistic.  This work closes the moment gap between the previously known finite $(3+\delta)$ moment condition for Hotelling moderate deviations and the finite third moment boundary.  Our proof adopts a geometric approach. Truncation at the finite-third-moment scale yields a truncated covariance matrix whose concentration error is negligible relative to the target \(\chi_d\) tail.  Uniform projected self-normalized moderate deviations are then averaged over directions and, after a whitening map with a ridge term, become a spherical cap transform of the radial distribution.  A local Abel inversion recovers the radial tail, and deterministic comparison removes the ridge and the truncation.
\end{abstract}

\maketitle

\section{Introduction}

Self-normalized statistics have been studied extensively in limit theory.  Their normal and moderate deviation approximations often remain valid under moment assumptions weaker than those in the classical Cram\'er theory for ordinary sums, because the observations that enlarge the numerator also enter the normalizing denominator.  This compensation is already present in Student's statistic and in the self-normalized sum.

The univariate theory is by now well developed.  Berry--Esseen bounds for Studentized statistics were obtained by Bentkus and G\"otze \cite{BentkusGotze1996}, and Gin\'e, G\"otze and Mason \cite{GineGotzeMason1997} characterized the asymptotic normality of Student's statistic.  Shao \cite{Shao1999} proved a Cram\'er type large deviation theorem for Student's statistic and identified the finite third moment boundary for the range $x=o(n^{1/6})$.  Jing, Shao and Wang \cite{JingShaoWang2003} developed a general Cram\'er theory for self-normalized sums of independent variables under finite $(2+\delta)$ moment assumptions.  Systematic accounts may be found in de la Pe\~na, Lai and Shao \cite{DeLaPenaLaiShao2009} and Shao \cite{Shao2013Survey}.  Related universal and refined moderate deviation results include Jing, Shao and Zhou \cite{JingShaoZhou2008} and Gao, Shao and Shi \cite{GaoShaoShi2022}.

For multivariate data, the analogue is Hotelling's \(T^2\) statistic, introduced by Hotelling \cite{Hotelling1931}. Its classical distribution theory is given in Muirhead \cite{Muirhead1982}.
For large and moderate deviations of Hotelling’s \(T^2\), Dembo and Shao \cite{DemboShao2006} obtained logarithmic estimates, and Liu and Shao \cite{LiuShao2013} established a Cramér-type relative-error result under a finite \((3+\delta)\)-moment condition.
This condition is stronger than the finite third moment condition sufficient in the scalar self-normalized theory. Whether the finite third moment condition is sufficient in the multivariate setting was therefore left open.
  The purpose of the present paper is to close this gap in fixed dimension by showing that the finite third moment already suffices for the full relative approximation in the range $o(n^{1/6})$.

Let
\[
 S_n=\sum_{j=1}^nX_j,
 \qquad A_n=\sum_{j=1}^nX_jX_j^\top,
 \qquad C_n=n^{-1}A_n,
 \qquad s_n=n^{-1/2}S_n.
\]
Our main object is  to study the Moore--Penrose self-normalized radius
\[
 U_{n,+}^2=s_n^\top C_n^+s_n=S_n^\top A_n^+S_n.
\]
So $U_{n,+}$ is the usual quadratic self-normalized radius when $C_n$ is nonsingular, and it is also well defined on singular samples.  We prove that, for every $a_n=o(n^{1/6})$,
\[
 \PP(U_{n,+}\ge x)=\PP(\chi_d\ge x)\{1+o(1)\}
\]
uniformly for $0\le x\le a_n$.  An exact Sherman--Morrison identity then transfers the result to the centered Hotelling statistic.

In what follows, we briefly outline the proof strategy.  The argument is not a reduction to an existing multivariate Hotelling theorem.  Instead, the proof starts with univariate self-normalized moderate deviations and extends them to the multivariate setting through a sequence of geometric arguments.  Specifically, Lemma~\ref{lem:scalar-jsw} provides the scalar estimate.  After truncation and ridge regularization, Lemma~\ref{lem:trunc-projected-tail} applies it uniformly to all projection directions.  Averaging those projected tails over directions gives the cap transform in Proposition~\ref{prop:trunc-cap-transform}.  To see why this transform contains the radial tail, fix a vector $y\in\R^d$.  The set of directions $u\in \mathbb{S}^{d-1}$ for which $u^\top y\ge t$ is a spherical cap.  In the moderate deviation window, Lemma~\ref{lem:cap-estimates} writes the threshold and radius as $x+s/x$ and $x+w/x$ and shows that the leading cap kernel is $(w-s)_+^{(d-1)/2}$, where $x_+:=\max\{0,x\}$ for all $x\in\R$.  This is the Abel kernel associated with integration over $(s,\infty)$, equivalently the kernel of a Riemann--Liouville fractional integral taken to the right.  Standard references include \cite{GorenfloVessella1991,SamkoKilbasMarichev1993}.  Lemmas~\ref{lem:abel-uniqueness} and~\ref{lem:abel-inversion}, together with the uniform form in Lemma~\ref{lem:uniform-abel-inversion}, then recover the radial tail from this local cap transform.

The main technical difficulty is empirical whitening at the finite third moment boundary.  Direct control of the surface Jacobian on rare events would require more than a third moment.  The truncation scheme in \Cref{sec:truncation} avoids this difficulty by using, for a prescribed sequence $x=x_n=o(n^{1/6})$, the level
\[
 \tau_n=\varepsilon_n\frac{\sqrt n}{1+x},
 \qquad
 \varepsilon_n=\left\{\frac{(1+x)^3}{\sqrt n}\right\}^{1/8}.
\]
Lemma~\ref{lem:trunc-bias} shows that this truncation has negligible bias and also gives
\[
 n\PP(\|X\|>\tau_n)
 \le \E\|X\|^3\frac{(1+x)^3}{\varepsilon_n^3\sqrt n}=o(1),
\]
so only a vanishing expected number of observations is discarded.  Lemma~\ref{lem:trunc-cov-concentration} shows that the truncated covariance concentrates with exponent $x^2/\varepsilon_n^2$, making the concentration error negligible relative to the target moderate deviation tail.  A ridge term $\alpha x^{-2}I_d$ then stabilizes the whitening map.  Proposition~\ref{prop:trunc-radial-tail} removes this ridge by deterministic comparison, and Lemma~\ref{lem:recursive-comparison}, based on the single vector deletion comparison in Lemma~\ref{lem:deletion}, restores the discarded observations.

The paper is organized as follows.  \Cref{sec:main-result} states the theorem and two deterministic facts.  \Cref{sec:truncation} introduces the truncated triangular array and proves a projected ridge tail estimate.  \Cref{sec:ridge-cap-transform} converts this projected estimate into a spherical cap transform with ridge regularization.  \Cref{sec:abel-inversion} proves the Abel inversion.  \Cref{sec:removing-ridge,sec:removing-truncation} remove the ridge and the truncation, respectively.  \Cref{sec:hotelling-transfer} transfers the self-normalized result to Hotelling's statistic.

\section{Main result}
\label{sec:main-result}

Throughout, $d\ge1$ is fixed and
\begin{equation}
\label{eq:basic}
 \E X=0,\qquad \E XX^\top=I_d,\qquad \kappa_3:=\E\|X\|^3<\infty.
\end{equation}
For a symmetric nonnegative matrix $A$, $A^+$ denotes its Moore--Penrose inverse.  The symbol $\|\cdot\|$ denotes Euclidean norm for vectors and spectral norm for matrices.  For symmetric matrices $A$ and $B$, the notation $A\preceq B$ means that $B-A$ is nonnegative definite.  Let $\Phi$ denote the standard normal distribution function and write $\bPhi=1-\Phi$.  Let $\chi_d$ denote the Euclidean norm of a standard normal vector in $\R^d$.  Throughout, $\mathbb{S}^{d-1}=\{u\in\R^d:\|u\|=1\}$ and $\sigma$ denotes normalized surface measure on $\mathbb{S}^{d-1}$.

Let
\[
 B_n=A_n-n^{-1}S_nS_n^\top,
 \qquad T_n^2=S_n^\top B_n^{-1}S_n
\]
when $B_n$ is nonsingular.  On the singular event, set $T_n=0$.

\begin{theorem}
\label{thm:main}
Assume \eqref{eq:basic}.  If $a_n=o(n^{1/6})$, then
\begin{equation}
\label{eq:main-U}
 \sup_{0\le x\le a_n}
 \left|\frac{\PP(U_{n,+}^2\ge x^2)}{\PP(\chi_d^2\ge x^2)}-1\right|\to0.
\end{equation}
The same conclusion holds for the centered Hotelling statistic, namely,
\begin{equation}
\label{eq:main-T}
 \sup_{0\le x\le a_n}
 \left|\frac{\PP(T_n^2\ge x^2)}{\PP(\chi_d^2\ge x^2)}-1\right|\to0.
\end{equation}
\end{theorem}

It is sufficient to prove the above result for an arbitrary sequence $x=x_n\to\infty$ with $x=o(n^{1/6})$.  The bounded range follows from the multivariate central limit theorem, $C_n\to I_d$ in probability, and continuity of the $\chi_d$ law.  The usual subsequence argument then gives uniformity.

 Before ending this section, we state two deterministic facts that will be used later.  The variational representation below rewrites the Moore--Penrose radius as the largest univariate self-normalized projection.  The same representation also leads to a single vector deletion comparison, which is the tool used in \Cref{sec:removing-truncation} to remove the truncation.

For deterministic vectors $z_1,\ldots,z_N\in\R^d$, let
\[
 \mathcal U(z_1,\ldots,z_N)
 =\left(\Bigl\{\sum_{i=1}^Nz_i\Bigr\}^\top
        \Bigl\{\sum_{i=1}^Nz_iz_i^\top\Bigr\}^+
        \Bigl\{\sum_{i=1}^Nz_i\Bigr\}\right)^{1/2}.
\]
Let $A=\sum_{i=1}^Nz_iz_i^\top$ and $s=\sum_{i=1}^Nz_i$.
Since $s\in\operatorname{range}(A)$, the generalized Cauchy--Schwarz
inequality
\[
   (u^\top s)^2\le (u^\top Au)(s^\top A^+s)
\]
holds for every $u$.  Equality is attained in the direction $A^+s$ when
$s\ne0$.  By homogeneity this gives
\begin{equation}
\label{eq:variational-U}
 \mathcal U(z_1,\ldots,z_N)
 =\sup_{\|u\|=1}
   \frac{|\sum_{i=1}^Nu^\top z_i|}
        {(\sum_{i=1}^N(u^\top z_i)^2)^{1/2}}
\end{equation}
with the convention $0/0=0$.

The variational formula also controls the effect of removing one vector.  The following single vector deletion comparison makes this precise and will be used in Lemma~\ref{lem:recursive-comparison} when the truncated observations are restored.

\begin{lemma}
\label{lem:deletion}
For $q>1$ and every $j\le N$,
\[
 \{\mathcal U(z_1,\ldots,z_N)\ge q\}
 \subset
 \{\mathcal U(z_1,\ldots,z_{j-1},z_{j+1},\ldots,z_N)
       \ge(q^2-1)^{1/2}\}.
\]
\end{lemma}

\begin{proof}
Fix a unit vector $u$ and write
\[
 s=\sum_{i\ne j}u^\top z_i,\qquad y=u^\top z_j,
 \qquad c=\sum_{i\ne j}(u^\top z_i)^2.
\]
If $c>0$, the elementary identity
\[
 \left(\frac{s^2}{c}+1\right)(c+y^2)-(s+y)^2
 =\left(c^{1/2}-\frac{sy}{c^{1/2}}\right)^2\ge0
\]
gives
\begin{equation}\label{syle}
 \frac{(s+y)^2}{c+y^2}\le \frac{s^2}{c}+1.
\end{equation}
If $c=0$, then $s=0$.  With the zero over zero convention, the ratio on the
left of \eqref{syle} is either zero or one, and the same inequality remains valid.  Taking
the supremum over $u$ and using \eqref{eq:variational-U} for the sample with
the $j$th vector deleted yields the deterministic inequality
\[
 \mathcal U(z_1,\ldots,z_N)^2
 \le
 \mathcal U(z_1,\ldots,z_{j-1},z_{j+1},\ldots,z_N)^2+1.
\]
The asserted event inclusion follows immediately.
\end{proof}

The deletion argument changes thresholds by small deterministic amounts, and the final Hotelling transfer compares nearby tail levels.  The following elementary estimates provide the needed control of the comparison tail.

\begin{lemma}
\label{lem:chi-tail}
For fixed $d\ge1$,
\begin{equation}
\label{eq:chi-tail-asymptotic}
 \PP(\chi_d\ge y)
 =\frac{2^{1-d/2}}{\Gamma(d/2)}
   y^{d-2}e^{-y^2/2}\{1+O(y^{-2})\},
 \qquad y\to\infty.
\end{equation}
Consequently, for all sufficiently large $x$ and all
$y\in[x/\sqrt2,x]$,
\begin{equation}
\label{eq:chi-tail-ratio-bound}
 \frac{\PP(\chi_d\ge y)}{\PP(\chi_d\ge x)}
 \le C_d\exp\{(x^2-y^2)/2\}.
\end{equation}
Moreover, if $h_x=o(x^{-1})$, then
\begin{equation}
\label{eq:chi-tail-local-stability}
 \frac{\PP(\chi_d\ge x+h_x)}{\PP(\chi_d\ge x)}\to1.
\end{equation}
\end{lemma}

\begin{proof}
The density of $\chi_d$ is
\[
 f_d(y)=\frac{2^{1-d/2}}{\Gamma(d/2)}
        y^{d-1}e^{-y^2/2},\qquad y>0.
\]
One integration by parts yields
\[
 \int_y^\infty t^{d-1}e^{-t^2/2}\,dt
 =y^{d-2}e^{-y^2/2}\{1+O(y^{-2})\}.
\]
This proves \eqref{eq:chi-tail-asymptotic}, uniformly when the argument is
bounded below by a quantity tending to infinity.  Applying it at $x$ and
$y\in[x/\sqrt2,x]$ proves \eqref{eq:chi-tail-ratio-bound}.  Finally,
\eqref{eq:chi-tail-asymptotic} and
\[
 (x+h_x)^2-x^2=2xh_x+h_x^2=o(1)
\]
give \eqref{eq:chi-tail-local-stability}.
\end{proof}

\begin{remark}
\label{rem:d-one}
When $d=1$, the assertion for $U_{n,+}$ is an immediate consequence of the
scalar self-normalized Cram\'er estimate.  Let $V_n^2=\sum_{i=1}^nX_i^2$ and
$R_n=S_n/V_n$, with $R_n=0$ on $\{V_n=0\}$.  Then
$U_{n,+}=|R_n|$.  If $x=x_n\to\infty$ and $x=o(n^{1/6})$, Lemma~\ref{lem:scalar-jsw}, applied to
$X_i$ and to $-X_i$, gives
\begin{eqnarray*}
 \PP(U_{n,+}\ge x)
 =\PP(R_n\ge x)+\PP(-R_n\ge x)
 &=&2\bPhi(x)\{1+o(1)\}\\
&=&\PP(\chi_1\ge x)\{1+o(1)\}.
\end{eqnarray*}
For bounded $x$, $V_n^2/n\to1$ in probability and
$S_n/\sqrt n\Rightarrow N(0,1)$, so $U_{n,+}\Rightarrow\chi_1$.  P\'olya's
theorem gives uniform absolute convergence of the distribution functions.  On
each compact  interval this is uniform relative convergence of the
tails because the limiting tail is bounded away from zero.  The usual
subsequence argument then proves \eqref{eq:main-U} for $d=1$.  The nonsingularity estimate and Sherman--Morrison transfer in \Cref{sec:hotelling-transfer} are valid
also in dimension one, and therefore give \eqref{eq:main-T} for $d=1$.
Consequently, the geometric proof below may and will assume $d\ge2$.
\end{remark}

\section{Truncation}
\label{sec:truncation}

From this point through \Cref{sec:removing-truncation} we work in the geometric case $d\ge2$.
Fix an arbitrary sequence $x=x_n\to\infty$ with $x=o(n^{1/6})$ and set
\begin{equation}
\label{eq:trunc-scales}
 h_n=\frac{(1+x)^3}{\sqrt n},\qquad
 \varepsilon_n=h_n^{1/8},\qquad
 \tau_n=\varepsilon_n\frac{\sqrt n}{1+x}.
\end{equation}
Then
\begin{equation}
\label{eq:scale-relations}
 \varepsilon_n\to0,\qquad \tau_n\to\infty,
 \qquad
 \frac{h_n}{\varepsilon_n^j}\to0\quad(j=1,2,3).
\end{equation}
Indeed, the last limits are immediate from $\varepsilon_n=h_n^{1/8}$, while
$\tau_n=n^{7/16}(1+x)^{-5/8}\to\infty$ because $1+x=o(n^{1/6})$.
Define
\[
 \widehat X_i=X_i\one\{\|X_i\|\le\tau_n\},
 \qquad \mu_n=\E\widehat X_i,
 \qquad M_n=\E(\widehat X_i\widehat X_i^\top).
\]
For $N\le n$, set
\[
 \widehat S_N=\sum_{i=1}^N\widehat X_i,\quad
 \widehat A_N=\sum_{i=1}^N\widehat X_i\widehat X_i^\top,\quad
 \widehat C_N=N^{-1}\widehat A_N,\quad
 \widehat s_N=N^{-1/2}\widehat S_N,
\]
and let $\widehat U_{N,+}=\mathcal U(\widehat X_1,\ldots,\widehat X_N)$.

The first step is to verify that this truncation has negligible bias and covariance error on the moderate deviation scale.  The same estimates also show that the number of deleted observations is small enough for the comparison in \Cref{sec:removing-truncation}.

\begin{lemma}
\label{lem:trunc-bias}
Uniformly over unit vectors $u$,
\begin{align}
 \|\mu_n\|&\le \kappa_3\tau_n^{-2},
 \label{eq:mu-bound}\\
 \|M_n-I_d\|&\le \kappa_3\tau_n^{-1},
 \label{eq:M-bound}\\
 \E|u^\top(\widehat X_i-\mu_n)|^3&\le C \kappa_3.
 \label{eq:third-centered-bound}
\end{align}
Moreover, let $\varrho_n=np_n=n\PP(\|X\|>\tau_n)$, then
\begin{equation}
\label{eq:rho-small}
 \varrho_n\le \kappa_3h_n\varepsilon_n^{-3}\to0.
\end{equation}
For $m_u=u^\top\mu_n$ and
\[
 \sigma_u^2
 =\operatorname{Var}(u^\top\widehat X_i)
 =u^\top M_nu-m_u^2,
\]
we have, uniformly for $\|u\|=1$, $N\in[n-x^2,n]$, and
$z\in[x/\sqrt2,2x]$,
\begin{equation}
\label{eq:bias-rates}
 z\sqrt N\,\|\mu_n\|\to0,
 \qquad
 z^2\|M_n-I_d\|\to0,
 \qquad
 z^2|\sigma_u^2-1|\to0.
\end{equation}
In particular, $1/2\le\sigma_u^2\le3/2$ for all sufficiently large $n$,
uniformly over $u$.
\end{lemma}

\begin{proof}
Since $\E X=0$,
\[
 \mu_n=-\E[X\one\{\|X\|>\tau_n\}],
\]
which gives \eqref{eq:mu-bound}.  Similarly,
\[
 I_d-M_n=\E[XX^\top\one\{\|X\|>\tau_n\}],
\]
and \eqref{eq:M-bound} follows.  The centered third moment bound follows from $|a-b|^3\le4(|a|^3+|b|^3)$ and Jensen's inequality.  Markov's inequality gives
\[
 np_n\le n\kappa_3\tau_n^{-3}=\kappa_3h_n\varepsilon_n^{-3},
\]
which proves \eqref{eq:rho-small}.  Finally,
\[
 z\sqrt N\,\|\mu_n\|
 \le C\frac{(1+x)^3}{\varepsilon_n^2\sqrt n}
 =C h_n\varepsilon_n^{-2}\to0,
\]
and
\[
 z^2\|M_n-I_d\|
 \le C\frac{(1+x)^3}{\varepsilon_n\sqrt n}
 =C h_n\varepsilon_n^{-1}\to0.
\]
Since
\[
 |\sigma_u^2-1|
 \le \|M_n-I_d\|+\|\mu_n\|^2
\]
and
\[
 z^2\|\mu_n\|^2
 \le \frac{\{z\sqrt N\,\|\mu_n\|\}^2}{N}\to0,
\]
the last assertion in \eqref{eq:bias-rates} follows from the first two
bounds.  The uniform lower and upper bounds for $\sigma_u^2$ are then
immediate.
\end{proof}

With the moment errors under control, the next step is concentration.  The following concentration lemma gives a probability bound far smaller than the target $\chi_d$ tail, so that later whitening errors can be ignored on the relative scale.

\begin{lemma}
\label{lem:trunc-cov-concentration}
For every fixed $\eta\in(0,1/2)$, uniformly for $N\in[n-x^2,n]$,
\begin{equation}
\label{eq:matrix-Bernstein}
 \PP(\|\widehat C_N-I_d\|>\eta)
 \le 2d\exp\left\{-c_\eta\frac{x^2}{\varepsilon_n^2}\right\}.
\end{equation}
If $Y_i=u^\top\widehat X_i$, $m_u=\E Y_i$, $\xi_i=Y_i-m_u$, and $\sigma_u^2=\E\xi_i^2$, then
\begin{equation}
\label{eq:scalar-Bernstein}
 \sup_{\|u\|=1}\PP\left(\left|N^{-1}\sum_{i=1}^N\xi_i^2-\sigma_u^2\right|>\eta\sigma_u^2\right)
 \le2\exp\left\{-c_\eta\frac{x^2}{\varepsilon_n^2}\right\}.
\end{equation}
In particular, both probabilities are $o\{\bPhi(z)\}$ uniformly for $z\in[x/\sqrt2,2x]$.
\end{lemma}

\begin{proof}
For all large $n$, \eqref{eq:M-bound} gives $\|M_n-I_d\|\le\eta/2$.
Moreover $M_n\preceq I_d$, because
$\widehat X_i\widehat X_i^\top\preceq X_iX_i^\top$ after taking expectations in quadratic forms.  The centered self adjoint summands
\[
 W_i=\widehat X_i\widehat X_i^\top-M_n
\]
satisfy $\|W_i\|\le2\tau_n^2$ for all sufficiently large $n$.
Using $(A-B)^2\preceq2A^2+2B^2$ for self-adjoint matrices and
\[
 (\widehat X_i\widehat X_i^\top)^2
 =\|\widehat X_i\|^2\widehat X_i\widehat X_i^\top,
\]
we obtain
\[
 \left\|\sum_{i=1}^N\E W_i^2\right\|
 \le C N\{\E\|\widehat X_i\|^4+\|M_n\|^2\}
 \le C N\tau_n \kappa_3.
\]
Here $\E\|\widehat X_i\|^4\le \tau_n\E\|X\|^3$ and the
bounded term $\|M_n\|^2$ is absorbed into $\tau_n\kappa_3$, since
$\tau_n\to\infty$ and $\kappa_3\ge(\E\|X\|^2)^{3/2}>0$.
The matrix Bernstein inequality \cite{Tropp2012} therefore gives
\[
 \PP(\|\widehat C_N-M_n\|>\eta/2)
 \le2d\exp\left\{-\frac{cN\eta^2}{\tau_n\kappa_3+\eta\tau_n^2}\right\}.
\]
Since $\tau_n\to\infty$, $N\ge n-x^2\sim n$, and
\[
 \frac{N}{\tau_n^2}
 =\frac{N}{n}\frac{(1+x)^2}{\varepsilon_n^2}
 \ge c\frac{x^2}{\varepsilon_n^2},
\]
the last formula proves \eqref{eq:matrix-Bernstein}.

For \eqref{eq:scalar-Bernstein}, we observe $|\xi_i|\le2\tau_n$ for all large $n$, and
\[
 \E\xi_i^4\le2\tau_n\E|\xi_i|^3\le C\tau_n \kappa_3.
\]
Moreover, $\sigma_u^2\in[1/2,3/2]$ uniformly by
Lemma~\ref{lem:trunc-bias}.  Thus the centered variables
$\xi_i^2-\sigma_u^2$ are bounded by $C\tau_n^2$ and have total variance at
most $CN\tau_n \kappa_3$. Applying Bernstein's inequality at level
$N\eta\sigma_u^2$ gives an exponent bounded below by
\[
 \frac{c_\eta N}{\tau_n\kappa_3+\tau_n^2}
 \ge c_\eta\frac{x^2}{\varepsilon_n^2},
\]
and hence proves \eqref{eq:scalar-Bernstein}, uniformly over $u$.
Finally, Mills' inequalities imply
\[
 \bPhi(z)=\exp\{-z^2/2+O(\log(1+z))\}.
\]
Since $z\le2x$ and $\varepsilon_n^{-2}\to\infty$, the exponential bounds above are $o\{\bPhi(z)\}$ uniformly over the stated range.
\end{proof}

The projected estimates reduce the multivariate problem to scalar self-normalized sums.  The next lemma states the univariate moderate deviation estimate that will be used immediately, uniformly over the projection direction $u$.

\begin{lemma}
\label{lem:scalar-jsw}
Let $\xi_1,\ldots,\xi_N$ be i.i.d.\ real random variables satisfying
$\E\xi_1=0$, $\E\xi_1^2=\sigma^2>0$ and
$\E|\xi_1|^3<\infty$.  Set
\[
 S_N=\sum_{i=1}^N\xi_i,
 \qquad V_N^2=\sum_{i=1}^N\xi_i^2,
 \qquad R_N=S_N/V_N,
\]
with $R_N=0$ on $\{V_N=0\}$, and define
\[
 L_N=\frac{\E|\xi_1|^3}{\sigma^3\sqrt N}.
\]
There are universal constants $c_0,C<\infty$ such that, whenever
\[
 (1+r)^3L_N\le c_0,
\]
\begin{equation}
\label{eq:scalar-jsw}
 \exp\{-C(1+r)^3L_N\}
 \le
 \frac{\PP(R_N\ge r)}{\bPhi(r)}
 \le
 \exp\{C(1+r)^3L_N\}.
\end{equation}
Consequently, for any triangular array whose standardized third absolute
moments are bounded uniformly over the array and over any auxiliary
parameters, one has
\begin{equation}
\label{eq:scalar-jsw-additive}
 \sup_{0\le r\le b_N}
 \left|\frac{\PP(R_N\ge r)}{\bPhi(r)}-1\right|\to0
 \qquad\text{for every } b_N=o(N^{1/6}).
\end{equation}
In particular, if $(1+b_N)^3L_N\to0$, then the supremum in
\eqref{eq:scalar-jsw-additive} is $O\{(1+b_N)^3L_N\}$.
\end{lemma}

\begin{proof}
This is the finite third moment specialization of the self-normalized
Cram\'er type moderate deviation theorem of Jing, Shao and Wang
\cite[Theorem~1.1]{JingShaoWang2003}, where the truncation quantity is at most $2(1+r)^3L_N$.
Indeed, the large summand contribution is
bounded by replacing $\xi_i^2$ with
$(1+r)|\xi_i|^3/(\sigma\sqrt N)$, and the truncated third moment contribution
is bounded directly by $\E|\xi_i|^3$.  Also $L_N\ge N^{-1/2}$ by Lyapunov's
inequality.  Thus, after choosing $c_0$ sufficiently small, the condition
$(1+r)^3L_N\le c_0$ verifies the range conditions in that theorem and gives
\eqref{eq:scalar-jsw}.  The statistic $S_N/V_N$ is invariant under
multiplication of all summands by the positive constant $1/\sigma$.

If the standardized third absolute moments are bounded, then
$L_N=O(N^{-1/2})$.  Hence, for every $b_N=o(N^{1/6})$,
\[
 \Delta_N:=\sup_{0\le r\le b_N}(1+r)^3L_N\to0.
\]
For all large $N$, the whole interval $0\le r\le b_N$ lies in the range of
\eqref{eq:scalar-jsw}.  Therefore \eqref{eq:scalar-jsw} and
$\max\{e^{C\Delta_N}-1,1-e^{-C\Delta_N}\}\le 2C\Delta_N$ give
\eqref{eq:scalar-jsw-additive} and the quantitative statement.
\end{proof}

We now apply the scalar estimate uniformly to univariate projections of the truncated array.  This step links the scalar theorem to the geometric cap transform by giving the tail of every projected ridge statistic at the same moderate deviation level.

\begin{lemma}
\label{lem:trunc-projected-tail}
Fix $\alpha\ge0$ and $L<\infty$.  Uniformly for
\[
 N\in[n-x^2,n],\qquad z\in[x/\sqrt2,2x],\qquad |s|\le L,
 \qquad u\in \mathbb{S}^{d-1},
\]
 let $t=z+s/z$ and
\[
 \widehat q_{N,\alpha,z}(u)
 =\frac{u^\top\widehat s_N}
 {(u^\top\widehat C_Nu+\alpha/z^2)^{1/2}},
\]
with a zero over zero ratio interpreted as zero when $\alpha=0$.
Then
\begin{equation}
\label{eq:trunc-projected-tail}
 \PP\{\widehat q_{N,\alpha,z}(u)\ge t\}
 =\bPhi((t^2+\alpha)^{1/2})\{1+o(1)\}.
\end{equation}
\end{lemma}

\begin{proof}
Use the notations in Lemma~\ref{lem:trunc-cov-concentration}, and let
\[
 S=\sum_{i=1}^N\xi_i,\qquad V^2=\sum_{i=1}^N\xi_i^2,
 \qquad
 R=\begin{cases}S/V,&V>0,\\0,&V=0.\end{cases}
\]
By Lemma~\ref{lem:trunc-bias}, the standardized variables
$\xi_i/\sigma_u$ are centered, have variance one, and have third absolute
moments bounded by a constant independent of $n$, $N$, and $u$.  Hence the
corresponding Lyapunov quantities $L_N$ in Lemma~\ref{lem:scalar-jsw} are
$O(N^{-1/2})$ uniformly.  Since $x=o(n^{1/6})$ and
$N\in[n-x^2,n]$, the interval $0\le r\le2x+K/x$ lies inside the
range $r\le c_0L_N^{-1/3}$ for all large $n$, and  the quantitative consequence
of Lemma~\ref{lem:scalar-jsw} gives, for every fixed $K<\infty$,
\begin{equation}
\label{eq:JSW-centered-trunc}
 \sup_{\substack{0\le r\le2x+K/x\\N\in[n-x^2,n]\\u\in \mathbb{S}^{d-1}}}
 \left|
 \frac{\PP(R\ge r)}{\bPhi(r)}-1
 \right|
 \le C_0\frac{(1+2x+K/x)^3}{\sqrt N}=o(1),
\end{equation}
where $C_0$ is independent of $n$, $N$, and $u$.

Fix $\eta\in(0,1/2)$ and set
\[
 E_{N,u}(\eta)=
 \left\{\left|\frac{V^2}{N\sigma_u^2}-1\right|\le\eta\right\}.
\]
On this event $V>0$.  Let $a=Nm_u/V$.  Define the deterministic quantity
\[
 \delta_n(\eta)=
 \sup_{\substack{N\in[n-x^2,n],\ z\in[x/\sqrt2,2x]\\u\in \mathbb{S}^{d-1}}}
 \frac{z\sqrt N\,|m_u|}{\sigma_u\sqrt{1-\eta}}.
\]
Lemma~\ref{lem:trunc-bias} implies
\begin{equation}
\label{eq:a-small}
 \delta_n(\eta)\to0,
 \qquad
 \sup_{\substack{z\in[x/\sqrt2,2x]\\u\in \mathbb{S}^{d-1}}}
 z^2|\sigma_u^2-1|\to0,
\end{equation}
and, on $E_{N,u}(\eta)$,
\begin{equation}
\label{eq:a-deterministic-bound}
 |a|\le\delta_n(\eta)/z.
\end{equation}
Writing $Y_i=\xi_i+m_u$ gives the exact identity
\begin{equation}
\label{eq:q-exact-centered}
 \widehat q_{N,\alpha,z}(u)
 =\frac{R+a}
 {\{1+2aR/N+a^2/N+N\alpha/(z^2V^2)\}^{1/2}}.
\end{equation}
Since $|R|\le\sqrt N$, there is a deterministic $\beta_n(\eta)\to0$,
uniform over the stated parameter set, such that on $E_{N,u}(\eta)$
\begin{equation}
\label{eq:denominator-small-bound}
 \left|\frac{2aR}{N}+\frac{a^2}{N}\right|
 \le \frac{\beta_n(\eta)}{z^2}.
\end{equation}
One may take
\[
 \beta_n(\eta)=
 \frac{4x\delta_n(\eta)}{\sqrt{n-x^2}}
 +\frac{\delta_n(\eta)^2}{n-x^2},
\]
which tends to zero because $x/\sqrt n\to0$.  Furthermore,
\[
 \gamma_-:=\frac{\alpha}{z^2\sigma_u^2(1+\eta)}
 \le\frac{N\alpha}{z^2V^2}
 \le\frac{\alpha}{z^2\sigma_u^2(1-\eta)}=:\gamma_+.
\]
For all sufficiently large $n$,
$1+\gamma_- -\beta_n(\eta)z^{-2}>0$.  Define the deterministic thresholds
\begin{align*}
 r_+&=t\left(1+\gamma_+
              +\frac{\beta_n(\eta)}{z^2}\right)^{1/2}
       +\frac{\delta_n(\eta)}{z},\\
 r_-&=t\left(1+\gamma_-
              -\frac{\beta_n(\eta)}{z^2}\right)^{1/2}
       -\frac{\delta_n(\eta)}{z}.
\end{align*}
Then \eqref{eq:q-exact-centered},
\eqref{eq:a-deterministic-bound}, and
\eqref{eq:denominator-small-bound} imply, on $E_{N,u}(\eta)$,
\[
 \{R\ge r_+\}
 \subset\{\widehat q_{N,\alpha,z}(u)\ge t\}
 \subset\{R\ge r_-\}.
\]
Consequently,
\[
 \PP(R\ge r_+)-\PP(E_{N,u}(\eta)^c)
 \le \PP\{\widehat q_{N,\alpha,z}(u)\ge t\}
 \le \PP(R\ge r_-)+\PP(E_{N,u}(\eta)^c).
\]
The scalar Bernstein bound makes the two error probabilities
$o\{\bPhi((t^2+\alpha)^{1/2})\}$ uniformly.  Moreover,
\eqref{eq:a-small} and the definition of $\beta_n(\eta)$ give
\[
 r_+^2=t^2+\frac{\alpha}{1-\eta}+o(1),
 \qquad
 r_-^2=t^2+\frac{\alpha}{1+\eta}+o(1),
\]
uniformly.  Here the factor $\sigma_u^{-2}=1+o(z^{-2})$ has been absorbed
into the $o(1)$ term.  Also $0<r_-\le r_+\le2x+C_\eta/x$ for all sufficiently
large $n$.  Thus \eqref{eq:JSW-centered-trunc} applies at the fixed
thresholds $r_\pm$.  Mills' ratio yields, uniformly over the parameter set,
\[
 \frac{\bPhi(r_+)}{\bPhi((t^2+\alpha)^{1/2})}
 \ge \exp\{-C_\alpha\eta+o(1)\},
 \qquad
 \frac{\bPhi(r_-)}{\bPhi((t^2+\alpha)^{1/2})}
 \le \exp\{C_\alpha\eta+o(1)\}.
\]
Letting first $n\to\infty$ and then
$\eta\downarrow0$ proves \eqref{eq:trunc-projected-tail}.
\end{proof}

The projected estimate becomes useful after averaging over directions.  This is the cap transform considered next.

\section{Ridge cap transform}
\label{sec:ridge-cap-transform}

Fix $\alpha>0$.  For the parameters of Lemma~\ref{lem:trunc-projected-tail}, define
\[
 \widehat C_{N,\alpha,z}=\widehat C_N+\alpha z^{-2}I_d,
 \qquad
 \widehat U_{N,\alpha,z}^2
 =\widehat s_N^\top\widehat C_{N,\alpha,z}^{-1}\widehat s_N.
\]
The Cauchy--Schwarz inequality in the inner product generated by
$\widehat C_{N,\alpha,z}$ gives
\[
 \widehat U_{N,\alpha,z}
 =\sup_{u\in \mathbb{S}^{d-1}}\widehat q_{N,\alpha,z}(u).
\]
Set
\[
 \widehat A_{t,\alpha,z}
 =\{u\in \mathbb{S}^{d-1}:\widehat q_{N,\alpha,z}(u)\ge t\}.
\]
Integrating \eqref{eq:trunc-projected-tail} over $u$ gives
\begin{equation}
\label{eq:mean-trunc-ellipsoid-cap}
 \E\sigma(\widehat A_{t,\alpha,z})
 =\bPhi((t^2+\alpha)^{1/2})\{1+o(1)\}.
\end{equation}
For $t\ge0$ and $r\ge0$, define
\begin{equation}
\label{eq:def-Hd}
 H_d(t,r)=\sigma\{v\in \mathbb{S}^{d-1}:rv^\top e_1\ge t\},
\end{equation}
where $e_1$ is any fixed unit vector.  Rotational invariance makes the
definition independent of the particular choice of $e_1$.

To turn the averaged projected tail into a statement about the radial statistic, we need to know how whitening changes surface measure on the sphere.  The next lemma gives that angular Jacobian and a Lipschitz bound near the identity.

\begin{lemma}
\label{lem:angular-jacobian}
Let $A$ be positive definite and define
\[
 F_A(u)=\frac{A^{1/2}u}{(u^\top Au)^{1/2}},
 \qquad u\in \mathbb{S}^{d-1}.
\]
Then the pushforward of the measure $\sigma$ under $F_A$ has density
\begin{equation}
\label{eq:angular-density}
 J_A(v)=\det(A)^{-1/2}(v^\top A^{-1}v)^{-d/2}
\end{equation}
with respect to $\sigma$.  Equivalently, for every Borel set
$D\subset \mathbb{S}^{d-1}$,
\[
 \sigma(F_A^{-1}(D))=\int_DJ_A(v)\,\sigma(dv).
\]
Moreover, $J_{cA}=J_A$ for every $c>0$, and for every fixed
$\eta_0<1$ there is $C_{d,\eta_0}<\infty$ such that
\begin{equation}
\label{eq:jacobian-Lipschitz}
 \|A-I_d\|\le\eta\le\eta_0
 \quad\Longrightarrow\quad
 \sup_{v\in \mathbb{S}^{d-1}}|J_A(v)-1|\le C_{d,\eta_0}\eta.
\end{equation}
\end{lemma}

\begin{proof}
Let $G\sim N(0,I_d)$ and $U=G/\|G\|$.  Then $U$ is independent of
$\|G\|$ and has law $\sigma$.  The direction of $A^{1/2}G$ is
$F_A(U)$.  Since $A^{1/2}G$ has density
\[
 (2\pi)^{-d/2}\det(A)^{-1/2}
 \exp\{-y^\top A^{-1}y/2\},
\]
polar integration $y=rv$ gives the angular density
\[
 \det(A)^{-1/2}(v^\top A^{-1}v)^{-d/2},
\]
the radial normalizing constant being the same as for $A=I_d$.  This
proves \eqref{eq:angular-density}.  The identity $J_{cA}=J_A$ follows
directly.  Finally, if $\|A-I_d\|\le\eta_0<1$, then the eigenvalues of
$A$ stay in a fixed compact subset of $(0,\infty)$.  The functions
$A\mapsto\log\det A$ and
$A\mapsto\log(v^\top A^{-1}v)$ are uniformly Lipschitz on this set,
uniformly in $v$.  Hence $|\log J_A(v)|\le C\eta$, and
\eqref{eq:jacobian-Lipschitz} follows.
\end{proof}

We apply Lemma~\ref{lem:angular-jacobian} with $A=\widehat C_{N,\alpha,z}$.  The whitening map
\[
 v=\frac{\widehat C_{N,\alpha,z}^{1/2}u}
 {(u^\top\widehat C_{N,\alpha,z}u)^{1/2}}
\]
transforms $\widehat A_{t,\alpha,z}$ into
\[
 \widetilde A_{t,\alpha,z}
 =\{v\in \mathbb{S}^{d-1}:
     v^\top\widehat C_{N,\alpha,z}^{-1/2}\widehat s_N\ge t\}.
\]
Writing
\[
 \widehat J_{N,\alpha,z}(v)
 =\det(\widehat C_{N,\alpha,z})^{-1/2}
   \{v^\top\widehat C_{N,\alpha,z}^{-1}v\}^{-d/2},
\]
we obtain
\begin{equation}
\label{eq:ellipsoid-cap-jacobian}
 \sigma(\widehat A_{t,\alpha,z})
 =\int_{\widetilde A_{t,\alpha,z}}
   \widehat J_{N,\alpha,z}(v)\,\sigma(dv).
\end{equation}
Since
\[
 \left\|\widehat C_{N,\alpha,z}^{-1/2}\widehat s_N\right\|
 =\widehat U_{N,\alpha,z},
\]
rotational invariance and \eqref{eq:def-Hd} give
\[
 \sigma(\widetilde A_{t,\alpha,z})
 =H_d(t,\widehat U_{N,\alpha,z}).
\]

The identity above shows that only the angular density separates the averaged ellipsoidal cap from an ordinary spherical cap.  The figure illustrates this whitening step, and the following proposition shows that the density error is negligible.

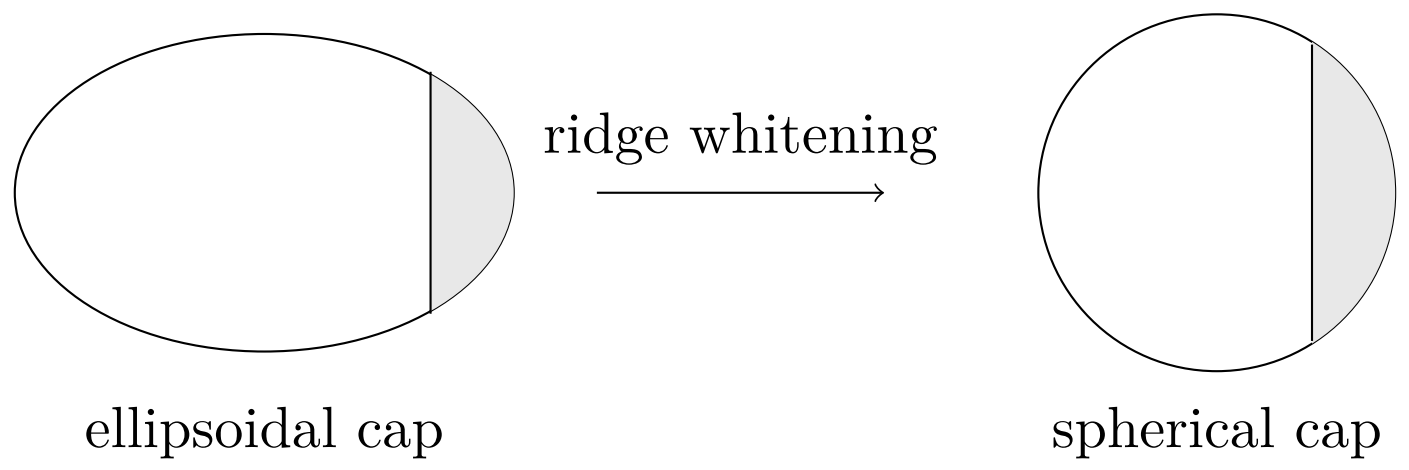
\begin{figure}[htbp]
\begin{center}
\begin{tikzpicture}[scale=1]
    \draw (-3.15,0) ellipse (1.65 and 1.05);
    \begin{scope}
        \clip (-3.15,0) ellipse (1.65 and 1.05);
        \fill[gray!18] (-2.05,-1.25) rectangle (-1.0,1.25);
    \end{scope}
    \draw (-2.05,-0.8) -- (-2.05,0.8);
    \node[below] at (-3.15,-1.28) {ellipsoidal cap};
    \draw[->] (-0.95,0) -- (0.95,0);
    \node[above] at (0,0.05) {ridge whitening};
    \draw (3.15,0) circle (1.18);
    \begin{scope}
        \clip (3.15,0) circle (1.18);
        \fill[gray!18] (3.78,-1.30) rectangle (4.45,1.30);
    \end{scope}
    \draw (3.78,-0.98) -- (3.78,0.98);
    \node[below] at (3.15,-1.28) {spherical cap};
\end{tikzpicture}
\caption{Ridge whitening converts an ellipsoidal excursion set into an ordinary spherical cap.  For the truncated array, the associated surface distortion is negligible on the moderate deviation scale.}
\label{fig:trunc-whitening}
\end{center}
\end{figure}

\begin{proposition}
\label{prop:trunc-cap-transform}
Under the uniform parameter range of Lemma~\ref{lem:trunc-projected-tail},
\begin{equation}
\label{eq:trunc-cap-transform}
 \E H_d(t,\widehat U_{N,\alpha,z})
 =\bPhi((t^2+\alpha)^{1/2})\{1+o(1)\}.
\end{equation}
\end{proposition}

\begin{proof}
Fix $\eta>0$  small enough that $\eta<1/2$ and
$C_d\eta<1/2$, and let
\[
 G_{N,\eta}=\{\|\widehat C_N-I_d\|\le\eta\}.
\]
Restricted to $G_{N,\eta}$,
\[
 \widehat C_{N,\alpha,z}
 =(1+\alpha z^{-2})
 \left\{I_d+\frac{\widehat C_N-I_d}{1+\alpha z^{-2}}\right\}.
\]
The scalar factor does not change the angular density, so
Lemma~\ref{lem:angular-jacobian} gives
\begin{equation}
\label{eq:J-good-bound}
 \sup_{v\in \mathbb{S}^{d-1}}
 |\widehat J_{N,\alpha,z}(v)-1|\le C_d\eta.
\end{equation}
Combining \eqref{eq:ellipsoid-cap-jacobian} and
\eqref{eq:J-good-bound}, we obtain on $G_{N,\eta}$
\[
 |\sigma(\widehat A_{t,\alpha,z})
   -H_d(t,\widehat U_{N,\alpha,z})|
 \le C_d\eta\,H_d(t,\widehat U_{N,\alpha,z}).
\]
On $G_{N,\eta}^c$, both cap measures lie in $[0,1]$.  Hence
\[
 \left|
 \E\sigma(\widehat A_{t,\alpha,z})
 -\E H_d(t,\widehat U_{N,\alpha,z})
 \right|
 \le C_d\eta\,\E H_d(t,\widehat U_{N,\alpha,z})
      +2\PP(G_{N,\eta}^c).
\]
Lemma~\ref{lem:trunc-cov-concentration} and Mills' ratio imply
\[
 \PP(G_{N,\eta}^c)
 =o\{\bPhi((t^2+\alpha)^{1/2})\}
\]
uniformly.  Together with \eqref{eq:mean-trunc-ellipsoid-cap}, these estimates imply
\[
 (1-C_d\eta)\E H_d(t,\widehat U_{N,\alpha,z})
 \le \bPhi((t^2+\alpha)^{1/2})\{1+o(1)\},
\]
so the mean cap is $O\{\bPhi((t^2+\alpha)^{1/2})\}$ uniformly.  Substituting
this bound back into the same inequality yields
\[
 \E H_d(t,\widehat U_{N,\alpha,z})
 =\bPhi((t^2+\alpha)^{1/2})
   \{1+O_d(\eta)+o(1)\}.
\]
Letting $n\to\infty$ first and then $\eta\downarrow0$ proves the
proposition.
\end{proof}

It remains to invert the cap transform on the moderate deviation scale.

\section{Abel inversion}
\label{sec:abel-inversion}

Let
\[
    m=d-1,
    \qquad
    \gamma=m/2,
    \qquad
    c_m=\frac{2^\gamma |B_m|}{|\mathbb{S}^{d-1}|},
\]
where $|B_m|$ is the volume of the unit ball in $\R^m$.  The power kernel that appears below is the Abel kernel associated with integration over $(s,\infty)$, or equivalently the kernel of a Riemann--Liouville fractional integral taken to the right.  Standard references include \cite{GorenfloVessella1991,SamkoKilbasMarichev1993}.  In the present proof this kernel arises from the geometry of small spherical caps.  The cap transform depends on the pair $(t,U)$ only through the relative position of the threshold and the radius.  In the moderate deviation window we write
\[
    t=x+s/x,
    \qquad
    U=x+w/x.
\]
Then the cap measure is of order $x^{-m}$, and after multiplying by $x^m$ its leading term depends only on the overshoot difference $w-s$.  The first lemma in this section identifies the Abel kernel uniformly and also controls the part of the transform coming from large overshoots.

\begin{figure}[htbp]
\begin{center}
\begin{tikzpicture}[scale=0.8]
    \draw (0,0) circle (1.65);
    \draw (0,0) circle (2.15);
    \begin{scope}
        \clip (0,0) circle (2.15);
        \fill[gray!18] (1.65,-2.30) rectangle (2.45,2.30);
    \end{scope}
    \draw (1.65,-1.43) -- (1.65,1.43);
    \draw[->] (0,0) -- (1.65,0);
    \draw[->] (0,0) -- (2.15,0);
    \draw[<->] (1.65,0.28) -- (2.15,0.28);
\end{tikzpicture}
\caption{On the overshoot scale, the threshold is $x+s/x$ and the radius is $x+w/x$.  The cap depends to first order on $(w-s)/x$, which yields the Abel kernel $(w-s)_+^\gamma$.  }
\label{fig:abel-geometry}
\end{center}
\end{figure}
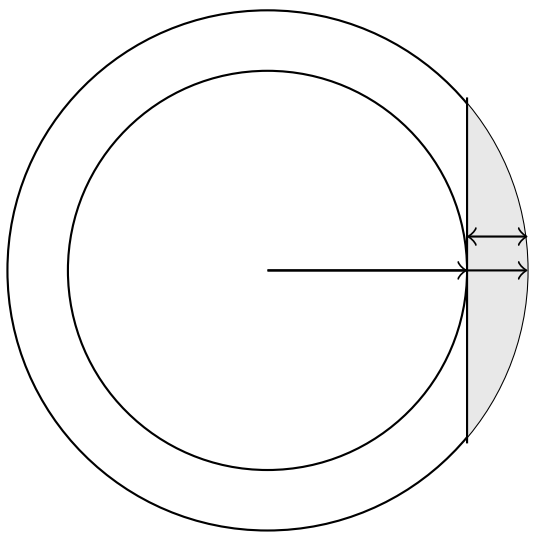

\begin{lemma}
\label{lem:cap-estimates}
For every compact set $K\subset\R^2$, as $x\to\infty$
\begin{equation}
\label{eq:small-cap-expansion}
    \sup_{(s,w)\in K}
    \left|
    \frac{x^m}{c_m}H_d(x+s/x,x+w/x)-(w-s)_+^\gamma
    \right|\to0.
\end{equation}
Moreover, for every fixed $L<\infty$ there exist constants
$R_0,C,c>0$ such that, for every fixed $R\ge R_0$, all sufficiently large
$x$ (depends on $R$), all $|s|\le L$, and all $w\ge R$,
\begin{equation}
\label{eq:cap-comparison-tail}
    H_d(x+s/x,x+w/x)
    \le C H_d(x+R/(2x),x+w/x),
\end{equation}
and
\begin{equation}
\label{eq:cap-tail-lower}
    H_d(x+(R-1)/x,x+w/x)
    \ge c x^{-m},
    \qquad w\ge R.
\end{equation}
The constants $R_0,C,c$ do not depend on $R$.
\end{lemma}

\begin{proof}
For $d\ge2$, we have  $\gamma=(d-1)/2>0$.  Let
$V$ have distribution $\sigma$ on $\mathbb{S}^{d-1}$ and take $e_1$ to be the
first coordinate vector.  For $0\le t<r$, rotational invariance and the
density of $V_1$ give
\[
 H_d(t,r)
 =\frac{|\mathbb{S}^{d-2}|}{|\mathbb{S}^{d-1}|}
   \int_{t/r}^1(1-y^2)^{(d-3)/2}\,dy.
\]
When $r>0$ and $t\ge r$, the cap has $\sigma$-measure zero.  Indeed, it is
empty if $t>r$, while for $t=r$ it consists of the single point $e_1$.

Set
\[
 t=x+s/x,\qquad r=x+w/x.
\]
For $(s,w)$ in a fixed compact set and all sufficiently large $x$, both
$t$ and $r$ are positive.  If $w\le s$, then $t\ge r$ and hence
$H_d(t,r)=0$.  Suppose now that $w>s$ and let
\[
 \Delta=1-\frac{t}{r}=\frac{w-s}{x^2+w}.
\]
Changing variables first by $y=1-q$ and then by $q=\Delta u$ yields the
exact representation
\[
 H_d(t,r)
 =\frac{|\mathbb{S}^{d-2}|}{|\mathbb{S}^{d-1}|}\Delta^\gamma
   \int_0^1u^{\gamma-1}(2-\Delta u)^{\gamma-1}\,du.
\]
For $0\le\Delta\le1$, the integral factor in this representation is finite
and depends continuously on $\Delta$.  At $\Delta=0$ its value is $2^{\gamma-1}/\gamma$.  Since
$|B_m|=|\mathbb{S}^{d-2}|/m$ and $m=2\gamma$,
\[
 \frac{|\mathbb{S}^{d-2}|}{|\mathbb{S}^{d-1}|}\frac{2^{\gamma-1}}{\gamma}
 =\frac{2^\gamma|B_m|}{|\mathbb{S}^{d-1}|}=c_m.
\]
If $K\subset\R^2$ is compact, then
\[
 \sup_{(s,w)\in K,\,w>s}\Delta\longrightarrow0
 \quad\text{and}\quad
 \sup_{(s,w)\in K}\left|\frac{x^2}{x^2+w}-1\right|\longrightarrow0.
\]
Consequently, for $w>s$,
\begin{align*}
 \frac{x^m}{c_m}H_d(x+s/x,x+w/x)
 &=(w-s)^\gamma
   \left(\frac{x^2}{x^2+w}\right)^\gamma \\
 &\quad\times
   \frac{|\mathbb{S}^{d-2}|}{c_m|\mathbb{S}^{d-1}|}
   \int_0^1u^{\gamma-1}(2-\Delta u)^{\gamma-1}\,du,
\end{align*}
where the last two terms converge uniformly to $1$ on $K$.  For
$w\le s$, both sides of the desired approximation are zero.  Since
$(w-s)_+^\gamma$ is bounded on $K$, this proves the absolute uniform
convergence in \eqref{eq:small-cap-expansion}.

We next derive bounds that remain valid when $w$ is unbounded.  For
$0\le t\le r$ with $r>0$, recall $\Delta=1-t/r\in[0,1]$.  The first change of variables
above gives
\[
 H_d(t,r)=\frac{|\mathbb{S}^{d-2}|}{|\mathbb{S}^{d-1}|}
 \int_0^\Delta q^{\gamma-1}(2-q)^{\gamma-1}\,dq.
\]
The continuous positive function $q\mapsto(2-q)^{\gamma-1}$ is bounded
above and below on $[0,1]$.  Hence there are constants
$0<c_1<C_1<\infty$, depending only on $d$, such that
\begin{equation}\label{eq:twosidesbound}
 c_1\Delta^\gamma\le H_d(t,r)\le C_1\Delta^\gamma,
 \qquad 0\le t\le r,\quad r>0.
\end{equation}
Fix $L<\infty$ and choose $R_0\ge4L+4$.  For $R\ge R_0$, $|s|\le L$
and $w\ge R$, define
\[
 \Delta_a(w)=1-\frac{x+a/x}{x+w/x}=\frac{w-a}{x^2+w}.
\]
For all sufficiently large $x$, the two thresholds corresponding to
$a=s$ and $a=R/2$ lie in $[0,x+w/x]$.  Moreover,
\[
 \frac{\Delta_s(w)}{\Delta_{R/2}(w)}
 =\frac{w-s}{w-R/2}
 \le4.
\]
Applying the two-sided bound (\ref{eq:twosidesbound}) to the two caps therefore gives
\begin{eqnarray*}
 H_d(x+s/x,x+w/x)
 \le C_1\Delta_s(w)^\gamma
 &\le& C_1 4^\gamma\Delta_{R/2}(w)^\gamma\\
&\le& \frac{C_1 4^\gamma}{c_1}
      H_d(x+R/(2x),x+w/x),
\end{eqnarray*}
which proves \eqref{eq:cap-comparison-tail} with a constant independent of
$R$.

Finally, fix $R\ge R_0$ and take $x$ sufficiently large that
$x\ge1$ and $x^2\ge\max\{R,2(R-1)\}$.  If $R\le w\le x^2$, then
\[
 \Delta_{R-1}(w)=\frac{w-R+1}{x^2+w}
 \ge\frac{1}{2x^2},
\]
and the same two-sided bound (\ref{eq:twosidesbound}) yields
\[
 H_d(x+(R-1)/x,x+w/x)
 \ge c_1 2^{-\gamma}x^{-2\gamma}
 =c_1 2^{-\gamma}x^{-m}.
\]
If $w>x^2$, then
\[
 \frac{x+(R-1)/x}{x+w/x}
 =\frac{x^2+R-1}{x^2+w}
 \le\frac34,
\]
so $\Delta_{R-1}(w)\ge1/4$ and
\begin{equation}\label{Hdx}
 H_d(x+(R-1)/x,x+w/x)
 \ge c_1 4^{-\gamma}
 \ge c_1 4^{-\gamma}x^{-m}.
\end{equation}
Thus \eqref{eq:cap-tail-lower} holds with
$c=c_1 4^{-\gamma}$.  The constants $R_0,C,c$ are independent of the
fixed value of $R$.
\end{proof}

The next lemma gives the uniqueness principle behind this inversion.  It will be used to identify the limiting overshoot measure from its Abel transform on the right.

\begin{lemma}
\label{lem:abel-uniqueness}
Let $\beta>0$ and let $\eta$ be a signed locally finite Borel measure on
$\R$.  Suppose that for some integer $k>\beta$ and every $A\in\R$,
\begin{equation}\label{intA}
    \int_A^\infty (1+w-A)^{k-1}|\eta|(dw)<\infty.
\end{equation}
If
\[
    I_-^\beta\eta(s)
:=\frac1{\Gamma(\beta)}
    \int_s^\infty (w-s)^{\beta-1}\eta(dw)=0,
    \qquad s\in\R,
\]
then $\eta=0$.
\end{lemma}

\begin{proof}
We observe
\[
 I_-^{k-\beta}(I_-^\beta\eta)(s)
 =\frac{1}{\Gamma(k-\beta)}
 \int_s^\infty
 (t-s)^{k-\beta-1}I_-^\beta\eta(t)\,dt.
\]
Substituting the definition of $I_-^\beta\eta(t)$ gives an integral over the
region $s\le t\le w$.  For each fixed $s$, the assumed
integrability makes the above double integral absolutely convergent, therefore, Fubini's theorem gives
\[
 I_-^{k-\beta}(I_-^\beta\eta)(s)
 =\frac{1}{\Gamma(k-\beta)\Gamma(\beta)}
 \int_s^\infty
 \left\{\int_s^w
 (t-s)^{k-\beta-1}(w-t)^{\beta-1}\,dt\right\}\eta(dw).
\]
With the substitution $t=s+(w-s)y$, the inner integral is
\[
 (w-s)^{k-1}\int_0^1y^{k-\beta-1}(1-y)^{\beta-1}\,dy
 =(w-s)^{k-1}
 \frac{\Gamma(k-\beta)\Gamma(\beta)}{\Gamma(k)}.
\]
Consequently,
\[
 I_-^{k-\beta}(I_-^\beta\eta)(s)
 =\frac1{\Gamma(k)}\int_s^\infty(w-s)^{k-1}\eta(dw)
 =I_-^k\eta(s).
\]
Thus $I_-^k\eta\equiv0$.  We now differentiate this integer-order
integral directly in the sense of distributions.  For fixed $w$, the
function
\[
 g_w(s)=\frac{(w-s)_+^{k-1}}{\Gamma(k)}
\]
satisfies $D_s^kg_w=(-1)^k\delta_w$.  If
$\varphi\in C_c^\infty(\R)$, the polynomial right-tail assumption permits
Fubini's theorem.  Since
$ I_-^k\eta(s)=\int_{\R}g_w(s)\,\eta(dw)$, the definition of a
distributional derivative gives
\begin{align*}
 \langle D^k I_-^k\eta,\varphi\rangle
 &=(-1)^k\int_{\R}
   \left\{\int_{\R}g_w(s)\,\eta(dw)\right\}
   \varphi^{(k)}(s)\,ds\\
 &=\int_{\R}
   \left\{(-1)^k\int_{\R}g_w(s)\varphi^{(k)}(s)\,ds\right\}
   \eta(dw)\\
 &=\int_{\R}\langle D_s^kg_w,\varphi\rangle\,\eta(dw)\\
 &=(-1)^k\int_{\R}\varphi(w)\,\eta(dw).
\end{align*}
Thus $D^kI_-^k\eta=(-1)^k\eta$.  Since $I_-^k\eta$ is identically zero,
we conclude that $\eta=0$.
\end{proof}

We can now pass from convergence of cap transforms to convergence of radial tails.

\begin{lemma}
\label{lem:abel-inversion}
Let $U=U_x\ge0$ be a sequence of random variables and let $x\to\infty$.
Suppose that $a>0$ and that, for every fixed $L<\infty$,
\begin{equation}
\label{eq:abel-hypothesis}
    \E H_d(x+s/x,U)
    =
    a e^{-s}\{1-\Phi(x)\}\{1+o(1)\}
\end{equation}
uniformly for $|s|\le L$.  Then, for every fixed $L<\infty$,
\[
    \PP(U\ge x+r/x)
    =
    a e^{-r}\PP(\chi_d\ge x)\{1+o(1)\}
\]
uniformly for $|r|\le L$.
\end{lemma}

\begin{proof}
Write
\[
 b_x=\frac{x^m\{1-\Phi(x)\}}{c_m},\qquad
 W_x=x(U-x),
\]
and define the finite positive measures
\[
 \nu_x(B)=\frac{\PP(W_x\in B)}{b_x},
 \qquad B\subset\R.
\]
For $w\ge-x^2$, set
\[
 h_{x,s}(w)=\frac{x^m}{c_m}H_d(x+s/x,x+w/x),
\]
and $h_{x,s}(w)=0$ for $w<-x^2$.  This convention outside
$[-x^2,\infty)$ is immaterial because $\nu_x$ is supported on that interval.
A direct change of variables with respect to the law of $W_x$ gives
\[
 \int h_{x,s}(w)\nu_x(dw)
 =\frac{\E H_d(x+s/x,U)}{1-\Phi(x)}.
\]
Hence the assumption (\ref{eq:abel-hypothesis}) is equivalent to
\begin{equation}
\label{eq:h-transform-hypothesis}
 \int h_{x,s}(w)\nu_x(dw)=a e^{-s}+o(1),
\end{equation}
locally uniformly in $s$.

We first obtain the compactness and right-tail bounds needed for the
inversion.  Fix $R\ge R_0$.  On the event $\{W_x\ge R\}$, equivalently
$\{U\ge x+R/x\}$, Lemma~\ref{lem:cap-estimates} gives, for all sufficiently
large $x$,
\[
 H_d(x+(R-1)/x,U)\ge c x^{-m}.
\]
Using the assumption (\ref{eq:abel-hypothesis}) at the fixed point $s=R-1$, we obtain
\[
 c x^{-m}\PP(W_x\ge R)
 \le \E H_d(x+(R-1)/x,U)
 \le 2a e^{1-R}\{1-\Phi(x)\}
\]
for all sufficiently large $x$.  Dividing by $b_x$ shows that, with a
constant $C$ independent of the fixed value of $R$,
\begin{equation}
\label{eq:nu-tail-tight}
 \nu_x([R,\infty))\le C e^{-R}.
\end{equation}
The index from which this estimate is valid may depend on $R$.

Let $[A,B]$ be compact and choose a fixed $s_0<A-1$.  By
\eqref{eq:small-cap-expansion},
\[
 \inf_{w\in[A,B]}h_{x,s_0}(w)
 \longrightarrow
 \inf_{w\in[A,B]}(w-s_0)^\gamma\ge1.
\]
Thus, for all sufficiently large $x$,
\[
 \nu_x([A,B])
 \le2\int h_{x,s_0}(w)\nu_x(dw)
 =2a e^{-s_0}+o(1).
\]
Each of the finitely many earlier $\nu_x$ is a finite measure, so enlarging
the constant covers those terms as well.  Therefore
\begin{equation}
\label{eq:nu-local-boundedness}
 \sup_x\nu_x([A,B])<\infty
\end{equation}
for every compact interval $[A,B]$.  Here and below, the supremum over $x$
are understood along the preassigned sequence $x\to\infty$ appearing
in the statement of lemma.

The bounds in \eqref{eq:nu-local-boundedness} imply relative compactness for
the vague topology.  To see this, for each integer $j\ge1$, we choose a
countable sup-norm dense set $\mathcal D_j\subset C_c((-j,j))$.  The union
$\mathcal D=\bigcup_j\mathcal D_j$ is countable, and for every
$f\in\mathcal D$ the numbers $\int f\,d\nu_x$ are bounded by
\eqref{eq:nu-local-boundedness}.  Starting from any subsequence and applying
a diagonal selection, we may pass to a further subsequence such that $\int f\,d\nu_x$
converges for every $f\in\mathcal D$.  Sup-norm approximation on the compact
support of an arbitrary $f\in C_c(\R)$, together with the corresponding
local mass bound, shows that $\int f\,d\nu_x$ also converges.  The resulting
functional on $C_c(\R)$ is positive and locally bounded, so the Riesz--Markov
theorem represents it as integration against a locally finite Borel measure
$\nu$.  Thus
\[
 \int f\,d\nu_x\longrightarrow\int f\,d\nu,
 \qquad f\in C_c(\R),
\]
and every subsequence has a further vaguely convergent subsequence.

We next show that overshoots far to the right make a uniformly negligible
contribution to the cap transform.  Fix $L<\infty$ and take
$R\ge R_0(L)$.  By \eqref{eq:cap-comparison-tail}, uniformly for
$|s|\le L$,
\begin{equation}
\label{eq:cap-tail-localization}
\begin{aligned}
 \int_{[R,\infty)}h_{x,s}(w)\nu_x(dw)
 &\le C\int_{[R,\infty)}h_{x,R/2}(w)\nu_x(dw)\\
 &\le C\int h_{x,R/2}(w)\nu_x(dw)\\
 &=C a e^{-R/2}+o(1).
\end{aligned}
\end{equation}
Here $R$ is fixed  before letting $x\to\infty$, and the final equality follows from
\eqref{eq:h-transform-hypothesis} at $s=R/2$.

Now let $\nu$ be a vague limit along an arbitrary subsequence.  Fix
$s\in\R$, let $L=|s|+1$, and choose a continuity point $R$ of $\nu$ with
$R>\max\{s,R_0(L)\}$.  Choose also a continuity point $b<s$.  Such choices
are possible because a locally finite measure has at most countably many
atoms.  For all sufficiently large $x$, $h_{x,s}(w)=0$ for $w\le s$ on the
support of $\nu_x$.  Hence
\[
 \int_{(-\infty,R]}h_{x,s}(w)\nu_x(dw)
 =\int_{[b,R]}h_{x,s}(w)\nu_x(dw).
\]
Vague convergence and the assumptions $\nu(\{b\})=\nu(\{R\})=0$ imply
weak convergence of the restricted measures on $[b,R]$.  Moreover,
\eqref{eq:small-cap-expansion} gives
\[
 \sup_{w\in[b,R]}
 \left|h_{x,s}(w)-(w-s)_+^\gamma\right|\longrightarrow0,
\]
and \eqref{eq:nu-local-boundedness} bounds
$\nu_x([b,R])$ uniformly.  Therefore
\begin{align*}
 &\left|
 \int_{[b,R]}h_{x,s}(w)\nu_x(dw)
 -\int_{[b,R]}(w-s)_+^\gamma\nu(dw)
 \right|\\
 &\quad\le
 \sup_{w\in[b,R]}|h_{x,s}(w)-(w-s)_+^\gamma|\,\nu_x([b,R])\\
 &\qquad+
 \left|
 \int_{[b,R]}(w-s)_+^\gamma\nu_x(dw)
 -\int_{[b,R]}(w-s)_+^\gamma\nu(dw)
 \right|
 \longrightarrow0.
\end{align*}
Combining this convergence with \eqref{eq:h-transform-hypothesis} and
\eqref{eq:cap-tail-localization} gives
\[
 0\le
 a e^{-s}-\int_s^R(w-s)^\gamma\nu(dw)
 \le C a e^{-R/2}.
\]
The left inequality follows because the omitted cap contribution over
$(R,\infty)$ is nonnegative. For the upper inequality it is bounded by the
integral over $[R,\infty)$ in \eqref{eq:cap-tail-localization}.

We claim that  (\ref{eq:nu-tail-tight}) also holds for $\nu$. In fact, if we let $R\ge R_0$ be
a continuity point of $\nu$, and choose continuity points $M>R$ tending to
infinity.  Vague convergence on the compact interval $[R,M]$ and
\eqref{eq:nu-tail-tight} give
\[
 \nu([R,M])
 =\lim_x\nu_x([R,M])
 \le\limsup_x\nu_x([R,\infty))
 \le C e^{-R}.
\]
Letting $M\to\infty$ and using monotone convergence yields
\begin{equation}
\label{eq:nu-limit-tail-bound}
 \nu([R,\infty))\le C e^{-R}
\end{equation}
for every sufficiently large continuity point $R$ of $\nu$.  In particular,
for every sufficiently large real $y$, one may choose a continuity point
$R_y\in[y-1,y]$ and obtain
\[
 \nu([y,\infty))\le\nu([R_y,\infty))
 \le C e^{-R_y}\le C e^{1-y}.
\]
Now let $R\to\infty$ through continuity points in the preceding
inequality.  The integrals on the left are
nondecreasing in $R$, and the error on the right tends to zero.  Hence
monotone convergence gives
\begin{equation}
\label{eq:limit-abel-transform}
 \int_s^\infty(w-s)^\gamma\nu(dw)=a e^{-s},
 \qquad s\in\R.
\end{equation}

Define
\[
 \nu_0(dw)=a\frac{e^{-w}}{\Gamma(1+\gamma)}\,dw.
\]
With the substitution $u=w-s$,
\[
 \int_s^\infty(w-s)^\gamma\nu_0(dw)
 =\frac{a e^{-s}}{\Gamma(1+\gamma)}
   \int_0^\infty u^\gamma e^{-u}\,du
 =a e^{-s},
\]
so $\nu_0$ has the same Abel transform as $\nu$.  Let
$\beta=1+\gamma$ and $\eta=\nu-\nu_0$.  Choose an integer $k>\beta$.
The bound \eqref{eq:nu-limit-tail-bound}, extended to all sufficiently large
real thresholds as above, implies for every $A\in\R$ that
\[
 \int_A^\infty(1+w-A)^{k-1}\nu(dw)<\infty.
\]
In fact,  splitting the integral into unit
intervals and using
$\nu([j,j+1])\le C e^{-j}$  prove the above assertion.  The same assertion is immediate for $\nu_0$.
Since $|\eta|\le\nu+\nu_0$, the polynomial right-tail condition in
Lemma~\ref{lem:abel-uniqueness} holds.
\eqref{eq:limit-abel-transform} and the corresponding identity for $\nu_0$
give
\[
 I_-^\beta\eta(s)
 =\frac1{\Gamma(\beta)}
   \int_s^\infty(w-s)^{\beta-1}\eta(dw)=0,
 \qquad s\in\R.
\]
Lemma~\ref{lem:abel-uniqueness} therefore yields $\eta=0$, and hence
$\nu=\nu_0$.

We have shown that every vaguely convergent subsequence has the same limit.
Together with the relative compactness proved above, this implies convergence
of the whole sequence
\begin{equation}
\label{eq:nu-limit-density}
 \nu_x\Rightarrow_v
 a\frac{e^{-w}}{\Gamma(1+\gamma)}\,dw.
\end{equation}

It remains to pass from vague convergence to convergence of right tails.
Fix $r\in\R$ and $\epsilon>0$.  By \eqref{eq:nu-tail-tight} and the explicit
right tail of $\nu_0$, choose $R>r$ so large that
\[
 \limsup_{x\to\infty}\nu_x([R,\infty))
 +\nu_0([R,\infty))<\epsilon.
\]
Because $\nu_0$ has no atoms, vague convergence gives
\[
 \nu_x([r,R])\longrightarrow\nu_0([r,R]).
\]
Since $[r,\infty)=[r,R]\mathbin{\dot\cup}(R,\infty)$, it follows that
\[
 \nu_x([r,\infty))\longrightarrow\nu_0([r,\infty))
 =\frac{a e^{-r}}{\Gamma(1+\gamma)}.
\]
This convergence is uniform for $r$ in compact intervals.  To see this
directly, fix $L<\infty$ and write
\[
 F_x(r)=\nu_x([r,\infty)),\qquad
 F(r)=\frac{a e^{-r}}{\Gamma(1+\gamma)}.
\]
Both $F_x$ and $F$ are nonincreasing, and $F$ is continuous.  Given
$\epsilon>0$, choose a partition
$-L=r_0<r_1<\cdots<r_N=L$ such that
$F(r_{j-1})-F(r_j)<\epsilon$ for every $j$.  Pointwise convergence at the
finitely many partition points and monotonicity then imply
\[
 \sup_{|r|\le L}|F_x(r)-F(r)|\le2\epsilon
\]
for all sufficiently large $x$.  Since
$\inf_{|r|\le L}F(r)>0$, the convergence is also uniform in relative error.
Consequently,
\begin{equation}
\label{eq:tail-from-nu}
 \PP(U\ge x+r/x)
 =\frac{x^m\{1-\Phi(x)\}}{c_m}
   \frac{a e^{-r}}{\Gamma(1+\gamma)}\{1+o(1)\}
\end{equation}
locally uniformly in $r$.

Finally, we identify the normalizing factor.  By
\eqref{eq:chi-tail-asymptotic},
\[
 \PP(\chi_d\ge x)
 =\frac{2^{1-d/2}}{\Gamma(d/2)}
   x^{d-2}e^{-x^2/2}\{1+o(1)\}.
\]
Furthermore,
\[
 |B_m|=\frac{\pi^{m/2}}{\Gamma(1+m/2)},
 \qquad
 |\mathbb{S}^{d-1}|=\frac{2\pi^{d/2}}{\Gamma(d/2)},
\]
so, using $m=2\gamma=d-1$,
\[
 c_m
 =\frac{2^{\gamma-1}}{\sqrt\pi}
   \frac{\Gamma(d/2)}{\Gamma(1+\gamma)}
 \quad\text{and therefore}\quad
 \frac{1}{\sqrt{2\pi}\,c_m\Gamma(1+\gamma)}
 =\frac{2^{1-d/2}}{\Gamma(d/2)}.
\]
Mills' ratio now gives
\begin{equation}
\label{eq:chi-tail-constant}
 \PP(\chi_d\ge x)
 =\frac{x^m\{1-\Phi(x)\}}
        {c_m\Gamma(1+\gamma)}\{1+o(1)\}.
\end{equation}
Substituting \eqref{eq:chi-tail-constant} into
\eqref{eq:tail-from-nu} proves the lemma.
\end{proof}

A compactness argument gives the uniform form needed in the sequel.

\begin{lemma}
\label{lem:uniform-abel-inversion}
Let $x_n\to\infty$ and let $\Theta_n$ be arbitrary index sets.  For each
$\theta\in \Theta_n$, let $z_{n,\theta}\in[x_n/\sqrt2,2x_n]$ and let
$U_{n,\theta}\ge0$ be a random variable.  Suppose that, for some $a>0$ and
for every fixed $L<\infty$,
\begin{equation}\label{eq:error}
 \sup_{\theta\in \Theta_n}\sup_{|s|\le L}
 \left|
 \frac{\E H_d(z_{n,\theta}+s/z_{n,\theta},U_{n,\theta})}{\bPhi(z_{n,\theta})}
 -ae^{-s}
 \right|\to0.
\end{equation}
Then, for every fixed $L<\infty$,
\[
 \sup_{\theta\in \Theta_n}\sup_{|r|\le L}
 \left|
 \frac{\PP(U_{n,\theta}\ge z_{n,\theta}+r/z_{n,\theta})}
 {ae^{-r}\PP(\chi_d\ge z_{n,\theta})}-1
 \right|\to0.
\]
\end{lemma}

\begin{proof}
If the conclusion failed, there would be a subsequence, indices
$\theta_n\in \Theta_n$, and bounded numbers $r_n$ for which the relative error in (\ref{eq:error}) stays away from zero.  Passing to a further subsequence, assume
$r_n\to r$.  Apply Lemma~\ref{lem:abel-inversion} to the single sequence
$z_n=z_{n,\theta_n}$ and $U_n=U_{n,\theta_n}$.  The assumption (\ref{eq:abel-hypothesis}) follows from
the assumed uniform cap transform convergence, and its conclusion contradicts
the selected failure.  This proves the uniform assertion.
\end{proof}

We return to the truncated statistic and remove the artificial ridge.

\section{Removing the ridge}
\label{sec:removing-ridge}

For $t=z+s/z$, Mills' ratio gives
\[
 \bPhi((t^2+\alpha)^{1/2})
 =e^{-s-\alpha/2}\bPhi(z)\{1+o(1)\}
\]
locally uniformly in $s$, uniformly for $z\in[x/\sqrt2,2x]$.  Combining
this with Proposition~\ref{prop:trunc-cap-transform} verifies the hypothesis
of Lemma~\ref{lem:uniform-abel-inversion}, with
$\theta=(N,z)$ and $a=e^{-\alpha/2}$.  Therefore, uniformly for fixed
$|r|\le L$,
\begin{equation}
\label{eq:trunc-ridge-tail}
 \PP(\widehat U_{N,\alpha,z}\ge z+r/z)
 =e^{-r-\alpha/2}\PP(\chi_d\ge z)\{1+o(1)\},
\end{equation}
for $N\in[n-x^2,n]$ and $z\in[x/\sqrt2,2x]$.

The estimate still contains the ridge parameter $\alpha$.  The next proposition removes it by comparing the ridge statistic with the Moore--Penrose statistic on the covariance concentration event.

\begin{proposition}
\label{prop:trunc-radial-tail}
Uniformly for $N\in[n-x^2,n]$ and $z\in[x/\sqrt2,2x]$,
\begin{equation}
\label{eq:trunc-radial-tail}
 \PP(\widehat U_{N,+}\ge z)
 =\PP(\chi_d\ge z)\{1+o(1)\}.
\end{equation}
\end{proposition}

\begin{proof}
Fix $\alpha>0$.  The vector $\widehat s_N$ belongs to the range of
$\widehat C_N$.  In a spectral decomposition of $\widehat C_N$, the
function $\lambda\mapsto(\lambda+\alpha z^{-2})^{-1}$ is bounded above by
$\lambda^{-1}$ on the positive eigenspace.  Hence
$\widehat U_{N,\alpha,z}\le\widehat U_{N,+}$.  Therefore
\eqref{eq:trunc-ridge-tail} with $r=0$ gives
\[
 \liminf
 \frac{\PP(\widehat U_{N,+}\ge z)}
      {\PP(\chi_d\ge z)}
 \ge e^{-\alpha/2},
\]
uniformly over $N\in[n-x^2,n]$ and $z\in[x/\sqrt2,2x]$.  Letting $\alpha\downarrow0$
gives the required lower bound.

For the upper bound, fix $\eta\in(0,1/2)$ and work on
$G_{N,\eta}=\{\|\widehat C_N-I_d\|\le\eta\}$.  On this event
$\widehat C_N\succeq(1-\eta)I_d$, and hence
\[
 \widehat C_N+\alpha z^{-2}I_d
 \preceq\left(1+\frac{\alpha}{(1-\eta)z^2}\right)\widehat C_N.
\]
Inverting the positive definite matrices reverses the order, so
\[
 \widehat U_{N,\alpha,z}
 \ge\left(1+\frac{\alpha}{(1-\eta)z^2}\right)^{-1/2}
      \widehat U_{N,+}.
\]
Consequently,
\begin{align*}
 \PP(\widehat U_{N,+}\ge z)
 &\le
 \PP\left\{\widehat U_{N,\alpha,z}\ge
 z\left(1+\frac{\alpha}{(1-\eta)z^2}\right)^{-1/2}\right\}
 +\PP(G_{N,\eta}^c)\\
 &=
 \PP\left\{\widehat U_{N,\alpha,z}\ge
 z-\frac{\alpha}{2(1-\eta)z}+O(z^{-3})\right\}
 +o\{\PP(\chi_d\ge z)\}.
\end{align*}
The last estimate uses Lemma~\ref{lem:trunc-cov-concentration}, which gives
$\PP(G_{N,\eta}^c)=o\{\bPhi(z)\}$ uniformly for $z\in[x/\sqrt2,2x]$.  By
\eqref{eq:chi-tail-asymptotic}, for fixed $d\ge2$ the normal tail is bounded by
$C_d\PP(\chi_d\ge z)$ in this range.  Hence
$\PP(G_{N,\eta}^c)=o\{\PP(\chi_d\ge z)\}$.  Applying
\eqref{eq:trunc-ridge-tail} with
$r=-\alpha/\{2(1-\eta)\}+O(z^{-2})$ yields
\[
 \limsup
 \frac{\PP(\widehat U_{N,+}\ge z)}
      {\PP(\chi_d\ge z)}
 \le
 \exp\left\{\frac{\alpha}{2(1-\eta)}-\frac\alpha2\right\}.
\]
Letting $\eta\downarrow0$ gives the upper bound one.  The argument is
uniform in $N,z$ because all preceding estimates are uniform.
\end{proof}

It remains to compare the truncated and original samples.

\section{Removing truncation}
\label{sec:removing-truncation}

Let
\[
 A_i=\{\|X_i\|>\tau_n\},\qquad p_n=\PP(A_i),\qquad \varrho_n=np_n.
\]
Recall from \eqref{eq:rho-small} that $\varrho_n\to0$,  $U_{N,+}=\mathcal U(X_1,\ldots,X_N)$ and $\widehat U_{N,+}=\mathcal U(\widehat X_1,\ldots,\widehat X_N)$.  In this section the hat always refers to the truncation level \eqref{eq:trunc-scales} chosen from the original pair $(n,x)$.

Lemma~\ref{lem:deletion}  allows us to restore the original observations.  Each large vector can lower the effective threshold by at most one unit in squared radius, and the next recursion keeps track of that loss.

\begin{lemma}
\label{lem:recursive-comparison}
Let $m=\lfloor x^2/2\rfloor$ and
\[
 (n)_\ell=n(n-1)\cdots(n-\ell+1),\qquad (n)_0=1.
\]
For all sufficiently large $n$,
\begin{align}
 \PP(U_{n,+}\ge x)
 &\le
 \sum_{\ell=0}^{m}(n)_\ell p_n^\ell
 \PP\left(\widehat U_{n-\ell,+}\ge\sqrt{x^2-\ell}\right)
 +(n)_{m+1}p_n^{m+1},
 \label{eq:recursive-upper}\\
 \PP(U_{n,+}\ge x)
 &\ge
 \PP(\widehat U_{n,+}\ge x)
 -np_n\,\PP(\widehat U_{n-1,+}\ge x).
 \label{eq:recursive-lower}
\end{align}
\end{lemma}

\begin{proof}
For $N\le n$ and $q>1$, let
\[
 F_N(q)=\PP\{\mathcal U(X_1,\ldots,X_N)\ge q\},
 \qquad
 \widehat F_N(q)=
 \PP\{\mathcal U(\widehat X_1,\ldots,\widehat X_N)\ge q\}.
\]
Splitting according to whether all $A_1^c,\ldots,A_N^c$ occur gives
\[
 F_N(q)
 \le \widehat F_N(q)
 +\sum_{j=1}^N
 \PP\{A_j,\mathcal U(X_1,\ldots,X_N)\ge q\}.
\]
On the event in the $j$th summand, Lemma~\ref{lem:deletion} implies
\[
 \mathcal U(X_1,\ldots,X_{j-1},X_{j+1},\ldots,X_N)
 \ge\sqrt{q^2-1}.
\]
The event involving the remaining observations is independent of $A_j$.
By exchangeability,
\begin{equation}
\label{eq:one-step-recursion}
 F_N(q)
 \le \widehat F_N(q)
 +Np_nF_{N-1}(\sqrt{q^2-1}).
\end{equation}
For $\ell\le m$, the successive thresholds
$\sqrt{x^2-\ell}$ are larger than one for all sufficiently large $n$.
Iterating \eqref{eq:one-step-recursion} $m+1$ times yields
\eqref{eq:recursive-upper}.  In the terminal term we simply bound the
remaining probability by one.

For the lower bound,
\[
 F_n(x)
 \ge
 \PP(\widehat U_{n,+}\ge x,A_1^c\cap\cdots\cap A_n^c).
\]
Hence
\[
 F_n(x)
 \ge\widehat F_n(x)
 -\sum_{j=1}^n\PP(A_j,\widehat U_{n,+}\ge x).
\]
On $A_j$ one has $\widehat X_j=0$, so the statistic
$\widehat U_{n,+}$ equals the truncated statistic formed from the other
$n-1$ observations.  Independence and exchangeability give
\[
 \PP(A_j,\widehat U_{n,+}\ge x)
 =p_n\,\PP(\widehat U_{n-1,+}\ge x),
\]
which proves \eqref{eq:recursive-lower}.
\end{proof}

We now combine the recursive comparison with the truncated radial estimate.

\begin{proof}[Completion of the proof of \eqref{eq:main-U} for $d\ge2$]
We first consider an arbitrary sequence $x=x_n\to\infty$ with
$x=o(n^{1/6})$, using the truncation level in
\eqref{eq:trunc-scales}.  Proposition~\ref{prop:trunc-radial-tail} gives
\begin{equation}
\label{eq:truncated-main-level}
 \PP(\widehat U_{n,+}\ge x)
 =\PP(\chi_d\ge x)\{1+o(1)\}.
\end{equation}
For $1\le\ell\le m=\lfloor x^2/2\rfloor$, let
$x_\ell=(x^2-\ell)^{1/2}$.  Then
$x/\sqrt2\le x_\ell\le x$, and
Proposition~\ref{prop:trunc-radial-tail} applies uniformly to
$\widehat U_{n-\ell,+}$ at threshold $x_\ell$.  In particular, for all
large $n$,
\[
 \PP(\widehat U_{n-\ell,+}\ge x_\ell)
 \le C\PP(\chi_d\ge x_\ell),
 \qquad 0\le\ell\le m.
\]
By Lemma~\ref{lem:chi-tail},
\begin{equation}
\label{eq:chi-deletion-ratio}
 \frac{\PP(\chi_d\ge x_\ell)}{\PP(\chi_d\ge x)}
 \le C_d e^{\ell/2},
 \qquad 1\le\ell\le m.
\end{equation}
Since $(n)_\ell p_n^\ell\le(np_n)^\ell=\varrho_n^\ell$ and
$\varrho_n\to0$, the terms with $1\le\ell\le m$ in
\eqref{eq:recursive-upper}, divided by $\PP(\chi_d\ge x)$, are bounded by
\[
 C\sum_{\ell=1}^{m}(e^{1/2}\varrho_n)^\ell=o(1).
\]
The  term $(n)_{m+1}p_n^{m+1}$ is negligible as well.  Indeed, for all sufficiently
large $n$, $\varrho_n\le e^{-2}$, and hence
\[
 (n)_{m+1}p_n^{m+1}
 \le\varrho_n^{m+1}
 \le e^{-2(m+1)}.
\]
Because $m=x^2/2+O(1)$ and
\eqref{eq:chi-tail-asymptotic} gives
$\log\PP(\chi_d\ge x)=-x^2/2+O(\log x)$,
\[
 \frac{(n)_{m+1}p_n^{m+1}}{\PP(\chi_d\ge x)}\to0.
\]
Combining these estimates with the $\ell=0$ term
\eqref{eq:truncated-main-level} in \eqref{eq:recursive-upper} gives
\begin{equation}
\label{eq:untruncated-upper}
 \limsup_{n\to\infty}
 \frac{\PP(U_{n,+}\ge x)}{\PP(\chi_d\ge x)}\le1.
\end{equation}

For the reverse inequality, Proposition~\ref{prop:trunc-radial-tail},
applied uniformly at $N=n$ and $N=n-1$, gives
\[
 \PP(\widehat U_{n-1,+}\ge x)
 \le C\PP(\chi_d\ge x).
\]
Therefore the error in \eqref{eq:recursive-lower} is at most
$C\varrho_n\PP(\chi_d\ge x)=o\{\PP(\chi_d\ge x)\}$.  Together with
\eqref{eq:truncated-main-level}, this yields
\begin{equation}
\label{eq:untruncated-lower}
 \liminf_{n\to\infty}
 \frac{\PP(U_{n,+}\ge x)}{\PP(\chi_d\ge x)}\ge1.
\end{equation}
Thus the desired relative approximation holds along every sequence
$x_n\to\infty$ with $x_n=o(n^{1/6})$.

It remains to justify both the bounded range and uniformity.  Since $d$ is
fixed and $\E\|X\|^2<\infty$, the entrywise law of large numbers gives
$C_n\to I_d$ in probability, while the multivariate central limit theorem
gives $s_n\Rightarrow Z\sim N(0,I_d)$.  Moreover $C_n^+\to I_d$ in
probability.  Indeed, on the event $\|C_n-I_d\|\le1/2$, the matrix $C_n$ is
invertible and
\[
 \|C_n^+-I_d\|=\|C_n^{-1}-I_d\|
 \le 2\|C_n-I_d\|.
\]
Thus Slutsky's theorem gives
\[
 U_{n,+}=(s_n^\top C_n^+s_n)^{1/2}\Rightarrow\|Z\|=\chi_d.
\]
Since the limiting distribution is continuous, P\'olya's theorem implies
uniform convergence of the distribution functions.  Consequently, for every
fixed $M<\infty$,
\[
 \sup_{0\le y\le M}
 \left|\frac{\PP(U_{n,+}\ge y)}{\PP(\chi_d\ge y)}-1\right|\to0,
\]
because $\inf_{0\le y\le M}\PP(\chi_d\ge y)>0$.

Finally, suppose the uniform assertion in \eqref{eq:main-U} failed for some
deterministic $a_n=o(n^{1/6})$.  There would be a sequence
$y_n\in[0,a_n]$ along which the relative error stays bounded away from zero.
A subsequence of $(y_n)$ is either bounded, contradicting the preceding
bounded range result, or tends to infinity, contradicting
\eqref{eq:untruncated-upper}--\eqref{eq:untruncated-lower}.  This proves
\eqref{eq:main-U}.
\end{proof}

We now prove the result for the centered Hotelling statistic.
\section{Hotelling transfer}
\label{sec:hotelling-transfer}

The identity covariance condition  provides the affine nondegeneracy needed for the transfer.

\begin{lemma}
\label{lem:affine-nondegeneracy}
Let $X$ be an $\R^d$-valued random vector with positive definite covariance.  There exists $\rho<1$ such that
\[
 \sup_H\PP(X\in H)\le\rho,
\]
where the supremum is over all affine hyperplanes $H\subset\R^d$.
\end{lemma}

\begin{proof}
Suppose the conclusion fails.  Then there are affine hyperplanes
\[
 H_k=\{y\in\R^d:u_k^\top y=b_k\},
 \qquad \|u_k\|=1,
\]
such that $\PP(X\in H_k)\to1$.  Choose $R<\infty$ with
$\PP(\|X\|\le R)>3/4$.  For all sufficiently large $k$, the events
$\{X\in H_k\}$ and $\{\|X\|\le R\}$ have an intersection of positive
probability.  Hence $H_k$ meets the closed ball of radius $R$, which implies
$|b_k|\le R$.   Using a subsequence argument, we may assume that
$u_k\to u$ and $b_k\to b$, where $\|u\|=1$.

Let $H=\{y:u^\top y=b\}$.  For every $y\notin H$, one has
$y\notin H_k$ for all sufficiently large $k$.  Fatou's lemma applied to the
complements therefore gives
\[
 \PP(X\notin H)
 \le \liminf_{k\to\infty}\PP(X\notin H_k)=0.
\]
Thus $u^\top X=b$ almost surely, so
$u^\top\operatorname{Cov}(X)u=0$.  This contradicts positive definiteness.
\end{proof}

We can now rule out singular centered covariance matrices on the tail scale.  The resulting exponential bound is stronger than every moderate deviation probability considered here.

\begin{lemma}
\label{lem:centered-singularity}
Under \eqref{eq:basic}, there are constants $c,C>0$ such that
\[
 \PP(B_n\text{ is singular})\le Ce^{-cn}.
\]
\end{lemma}

\begin{proof}
Since
\[
 B_n=\sum_{i=1}^n(X_i-\bar X_n)(X_i-\bar X_n)^\top,
 \qquad \bar X_n=S_n/n,
\]
the rank of $B_n$ equals the dimension of the affine span of
$X_1,\ldots,X_n$, which are supposed to be observed in order.  Conditional on
the first $j$ observations, if their affine span has dimension at most
$d-1$, it is contained in an affine hyperplane.  Independence and Lemma~\ref{lem:affine-nondegeneracy} imply that the conditional probability that
$X_{j+1}$ remains in the current affine span is at most $\rho$.

Full affine dimension $d$ requires at least $d$ dimension increasing steps
after $X_1$.  For any prescribed pattern containing exactly $k$ such steps,
the probability of the remaining $n-1-k$ failures is at most
$\rho^{n-1-k}$.  The probabilities of the increasing steps are bounded by
one.  Summing over patterns therefore gives
\[
 \PP(B_n\text{ is singular})
 \le\sum_{k=0}^{d-1}\binom{n-1}{k}\rho^{n-1-k}.
\]
If $\rho=0$, the last expression vanishes for all $n>d$.  If
$0<\rho<1$, it is bounded by
$C n^{d-1}\rho^{n-d}\le Ce^{-cn}$.  This proves the claim.
\end{proof}

On this high probability nonsingular event in Lemma \ref{lem:centered-singularity}, the two statistics are linked by an exact identity.

\begin{proof}[Completion of the proof of \eqref{eq:main-T}]
On the event that $B_n$ is nonsingular, it is positive definite and
\[
 A_n=B_n+n^{-1}S_nS_n^\top
\]
is positive definite as well.  Thus $U_{n,+}=U_n$, where
$U_n^2=S_n^\top A_n^{-1}S_n$.  Let $v_n=S_n/\sqrt n$.  Since
$B_n=A_n-v_nv_n^\top$, we have
\[
 A_n^{-1/2}B_nA_n^{-1/2}
 =I_d-w_nw_n^\top,
 \qquad w_n=A_n^{-1/2}v_n.
\]
The left side is positive definite, so
$1-\|w_n\|^2=1-v_n^\top A_n^{-1}v_n>0$.  The
Sherman--Morrison identity is therefore applicable and gives
\[
 B_n^{-1}
 =A_n^{-1}
  +\frac{A_n^{-1}v_nv_n^\top A_n^{-1}}
         {1-v_n^\top A_n^{-1}v_n}.
\]
As $v_n^\top A_n^{-1}v_n=U_n^2/n$, multiplication on both sides by $S_n$
yields the exact identity
\begin{equation}
\label{eq:SM}
 T_n^2=\frac{U_n^2}{1-U_n^2/n}.
\end{equation}
The map on the right is increasing for $0\le U_n^2<n$, and hence, on the
nonsingular event,
\begin{equation}
\label{eq:T-U-event-equivalence}
 \{T_n\ge x\}
 =\left\{U_{n,+}\ge y_n(x)\right\},
 \qquad
 y_n(x)=\frac{x}{\sqrt{1+x^2/n}}.
\end{equation}

Consider first an arbitrary sequence $x=x_n\to\infty$ with
$x=o(n^{1/6})$.  Then
\[
 y_n(x)
 =x-\frac{x^3}{2n}\{1+o(1)\}
 =x+o(x^{-1}),
\]
because $x^4/n\to0$.  In particular, $y_n(x)\to\infty$ and
$y_n(x)=o(n^{1/6})$.  The established assertion \eqref{eq:main-U} for $U_{n,+}$ and
Lemma~\ref{lem:chi-tail} therefore give
\[
 \PP\{U_{n,+}\ge y_n(x)\}
 =\PP\{\chi_d\ge y_n(x)\}\{1+o(1)\}
 =\PP(\chi_d\ge x)\{1+o(1)\}.
\]
Moreover, by Lemma~\ref{lem:centered-singularity},
\[
 \frac{\PP(B_n\text{ is singular})}{\PP(\chi_d\ge x)}\to0,
\]
since $x^2=o(n)$ and
$\log\PP(\chi_d\ge x)=-x^2/2+O(\log x)$.  Because
\eqref{eq:T-U-event-equivalence} holds on the nonsingular event,
\[
 \left|\PP(T_n\ge x)-\PP\{U_{n,+}\ge y_n(x)\}\right|
 \le \PP(B_n\text{ is singular}).
\]
The preceding formulas prove the relative approximation for every such
diverging sequence.

For bounded thresholds, note that
\[
 \frac{B_n}{n}=C_n-\frac{s_ns_n^\top}{n}\to I_d
 \quad\text{in probability},
\]
while $s_n\Rightarrow N(0,I_d)$.  On the event
$\|B_n/n-I_d\|\le1/2$, the matrix $B_n$ is nonsingular and
\[
 \|(B_n/n)^{-1}-I_d\|
 \le2\|B_n/n-I_d\|.
\]
The complement has probability tending to zero.  On this event,
\[
 \left|T_n^2-\|s_n\|^2\right|
 \le \|(B_n/n)^{-1}-I_d\|\,\|s_n\|^2.
\]
Since $s_n=O_{\mathbb P}(1)$, it follows that
$T_n=\|s_n\|+o_{\mathbb P}(1)$.  Slutsky's theorem gives
$T_n\Rightarrow\chi_d$.
Since the distribution of $\chi_d$ is continuous, P\'olya's theorem gives,
for every fixed $M<\infty$,
\[
 \sup_{0\le y\le M}
 \left|\PP(T_n\ge y)-\PP(\chi_d\ge y)\right|\to0.
\]
The tail probability $\PP(\chi_d\ge y)$ is bounded away from zero for
$0\le y\le M$.  Hence
\[
 \sup_{0\le y\le M}
 \left|
 \frac{\PP(T_n\ge y)}{\PP(\chi_d\ge y)}-1
 \right|\to0.
\]

Suppose that \eqref{eq:main-T} does not hold.  Then there is a sequence
$x_n\in[0,a_n]$ along which the relative error stays away from zero.
Along a subsequence, either $x_n$ is bounded or $x_n\to\infty$.  The first
case contradicts the preceding display, and the second contradicts the
result proved for $x_n=o(n^{1/6})$.  Therefore
\eqref{eq:main-T} holds.
\end{proof}

\end{document}